# Pascal's Triangles in Abelian and Hyperbolic Groups

*Michael Shapiro*

We are used to imagining Pascal's triangle as extending forever downwards from a vertex located at the top. But it is interesting to see it as occupying the first quadrant of the plane with it's vertex at $(0,0)$. Imagine further that the plane is made of graph paper — that is, that we have embedded into it the Cayley graph of $\mathbb{Z} \times \mathbb{Z}$ with respect to the standard generating set. If we place the entries of Pascal's triangle at the vertices of this Cayley graph, they now measure something about this graph. The entry at each point gives the number of geodesics from $(0,0)$ to that point.

This leads us to the following definition.

**Definition.** Suppose $\Gamma = \Gamma_\mathcal{G}(G)$ is the Cayley graph of $G$ with respect to the generating set $\mathcal{G}$. The **Pascal's function**, $p = p_\mathcal{G} : G \to \mathbb{Z}$ is given by

$$p(g) = \#\{\text{geodesics from } 1 \text{ to } g \text{ in } \Gamma_\mathcal{G}(G)\}.$$

This definition can be extended to any graph. We will only be interested in Cayley graphs. Conversations with several emminent geometric group theorists and combinatorists suggest that surpisingly little is known about these. I wish to thank Jim Cannon for his kind encouragement.

Let us be more specific about notation. We will take the generating set $\mathcal{G}$ to be a set which bijects to a subset of $G$ closed under inversion. The elements of $\mathcal{G}$ can be multiplied together in $\mathcal{G}^*$, the free monoid on $\mathcal{G}$ to form words. Their images can be multiplied together in $G$. The map taking words to their values in $G$ is a monoid homomorphism. We will denote it by $w \mapsto \overline{w}$. For a word $w$, $\ell(w)$ denotes its length. For a group element $g$, $\ell(g) = \ell_\mathcal{A}(g)$ denotes its length, i.e., the length of the shortest $\mathcal{A}$ word which evaluates to $g$. A word $w$ is geodesic if $\ell(w) = \ell(\overline{w})$.

Given $A \subset G$ and $\mathcal{G}$ generating $G$ we say that $A$ is **totally geodesic** if every $\mathcal{G}$ geodesic for an element of $A$ lies entirely in $A$. We will be interested in the case where $A$ is a subgroup or a submonoid. If $A$ is totally geodesic subgroup or submonoid, then $\mathcal{A} = \mathcal{G} \cap A$ is a generating set for $A$, and the following is immediate:

**Proposition.** *If $A < G$ is totally geodesic with respect to $\mathcal{G}$ then $p_\mathcal{A} = p_\mathcal{G}\big|_A$.* □

The Pascal's function $p_\mathcal{G}$ can depend very strongly on $\mathcal{G}$. For example, consider $\mathbb{Z}$. If we take the generating set consisting of a single generator, then $p$ is identically 1. However, if we take the generating set $\mathbb{Z} = \langle t, s \mid \overline{s} = \overline{t}^{10} \rangle$, then a number of the form $g = \bar{t}^{10k+5}$ with $k > 0$ has a Pascal's function which goes up rather quickly as a function of $k$. Specifically, if $n_j$ is the number of ordered partitions of 5 in $j$ blocks, $1 \leq j \leq 5$, then $p(g) = \sum_{j=1}^{5} \binom{k+1}{n_j}$. (This is because each geodesic for $g$ will consist of $k$ $s$'s together with 5 $t$'s distributed around and between them.) This goes up with the fifth power of $k$. Now in some sense this not too bad. In fact, given any generating set for $\mathbb{Z}$, there is a finite index subgroup (namely the one generated by the largest generator) which is totally geodesic, and the Pascal's function on this subgroup is identically 1. However, the dependence becomes more "ineradicable" if we turn to a free group of rank greater than 1. Once again if we take a basis, the Pascal's function is identically 1. Now consider $F_2$, the free group of rank two with the generating set $\langle x, y, a, b, c \mid \bar{a} = \bar{x}^3, \bar{b} = \bar{y}^3, \bar{c} = \overline{x^3 y^2} \rangle$. Then $\ell(\overline{x^3 y^3}) = 2$, and indeed, $\ell((\overline{x^3 y^3})^k) = 2k$, and $p((\overline{x^3 y^3})^k) = 2^k$. Any finite index subgroup must meet this subgroup, and thus the dependence on generating set will not go away by passing to a finite index subgroup.

## Abelian groups

**Proposition.** *Suppose that $G = A \times B$ and the $\mathcal{G} = \mathcal{A} \times \{1\} \cup \{1\} \times \mathcal{B}$, where $\mathcal{A}$ and $\mathcal{B}$ are generating sets for $A$ and $B$ respectively. Then*

$$p_{\mathcal{G}}(a,b) = \binom{\ell_{\mathcal{A}}(a) + \ell_{\mathcal{B}}(b)}{\ell_{\mathcal{A}}(a)} p_{\mathcal{A}}(a) p_{\mathcal{B}}(b).$$

*Proof.* A $\mathcal{G}$ geodesic $w$ for $(a,b)$ determines an $\mathcal{A}$ geodesic $w_a$ for $a$ and a $\mathcal{B}$ geodesic $w_b$ for $b$. Given the pair $w_a$ and $w_b$ there are exactly $\binom{\ell_{\mathcal{A}}(a) + \ell_{\mathcal{B}}(b)}{\ell_{\mathcal{A}}(a)}$ ways of combining them into a $\mathcal{G}$ geodesic for $(a,b)$. $\square$

This shows how to recover the standard Pascal's triangle from the Pascal's function for $\mathbb{Z}$ with respect to a single generator, or indeed how to find the Pascal's function of a finitely generated free abelian group with respect to a basis.

There is a sense in which the Pascal's function for $\mathbb{Z}^n$ with respect to a basis are the "prototype" Pascal's functions for abelian groups.

Let $A$ be an abelian group and let $\mathcal{A}$ be a generating set. We will say that a subset $S = \{a_{i_1}, \ldots, a_{i_k}\} \subset \mathcal{A}$ is **compatible** if for any $N$ there is a geodesic $w_N$ containing at least $N$ of each letter of $S$.

**Proposition.** *Let $S$ be a maximal compatible set, and let $M = M(S) = \overline{S^*}$ be the submonoid of $A$ generated by $\bar{S}$. Then $S^*$ is exactly the set of geodesics evaluating into $M$. Consequently the map*

$$S^* \to \mathbb{Z}_{\geq 0}^k \xrightarrow{\pi} M$$

takes $\mathbb{Z}_{\geq 0}^k$-geodesics to $M$-geodesics and for $a \in M$

$$p_{\mathcal{A}}(a) = \sum_{g \in \pi^{-1}(a)} p_{\mathbb{Z}^k}(g).$$

*Proof.* We first check that $S^*$ contains only geodesics. To see this, observe that the geodesics of an abelian group are closed under permutation and the geodesics of any group are closed under passing to subwords.

Next, we check that $S^*$ exhausts the geodesics of $M$. Suppose to the contrary that it does not. Then there is a geodesic $u$ evaluating into $M$ containing some letter (say, $a$) not in $S$. Let $v \in S^*$ be an $S$-geodesic with $\overline{u} = \overline{v}$, and let $w$ be any $S$ word containing all letters of $S$. Then for any $N$, $v^N w^N$ is a geodesic. But $\ell(u) = \ell(v)$ so $u^N w^N$ is also geodesic and contains at least $N$ instances of each letter of $S \cup \{a\}$. This contradicts the maximality of $S$. □

We can discover the maximal compatible sets via the use of translation lengths. For each element $g \in A$, we take the **translation length** $\tau(g)$ to be

$$\tau(g) = \lim_{j \to \infty} \frac{\ell(g^j)}{j}.$$

We consider the free abelian portion of $A$ to be $\mathbb{Z}^n \subset \mathbb{Q}^n \subset \mathbb{R}^n$. For each element of $A$, we have just defined the translation length. Given $q \in \mathbb{Q}^n$, there is $m$ so that $mq \in \mathbb{Z}^n$ and we define $\tau(q) = \frac{\tau(mq)}{m}$. This is independent of choice of $m$. Finally, we can extend $\tau$ to $\mathbb{R}^n$ by continuity. (For details, see [**NS**].) We take

$$C = \{x \in \mathbb{R}^n \mid \tau(x) \leq 1\}.$$

In the case where $A = \mathbb{Z}^n$, this is the convex hull of $\overline{\mathcal{A}} \subset \mathbb{R}^n$.

In the case where $A = \mathbb{Z}^n \times F$ with $F$ finite, $C$ is the convex hull of a related object. We take $f = \#F$ and let

$$W = \{w \in \mathcal{A}^* \mid \overline{w} \in \mathbb{Z}^n \text{ and } 0 < \ell(w) \leq f\}.$$

We take

$$V = \{\frac{\overline{w}}{\ell(w)} \mid w \in W\}.$$

If $a$ is a letter of $w \in W$, we say that $a$ **appears** at $\frac{\overline{w}}{\ell(w)} \in \mathbb{Q}^n$.

**Proposition.** *$C$ is the convex hull of $V$. $S$ is compatible if and only if all the elements of $S$ appear on a common face of $C$. $S$ is maximal compatible if and only if all the elements of $S$ appear on a maximal face of $C$.*

Notice that an element of $S$ may appear on the boundary of $C$ without being a vertex of $C$.

*Proof.* The first part of the proposition is a special case of [**NS**] Lemma 5.3. While [**NS**] are dealing with virtually abelian groups, our groups are abelian, so we can simplify the situation by taking

$$W' = \{a^f \mid a \in A\}$$

and
$$V' = \{\frac{1}{f}\bar{a}^f \mid a \in A\}$$
since in this case the convex hull of $V'$ is identical to the convex hull of $V$. Now each element of $S$ appears at exactly one point of $V'$.

We consider a set $S \subset \mathcal{A}$ and investigate when this is compatible. This fails to be compatible if and only if there is some word (which we write additively) $m_1 a_1 + \cdots + m_j a_j$ with each $a_i \in S$ which can be shortened, say as
$$m_1 a_1 + \cdots + m_j a_j = m'_1 a_1 + \cdots + m'_j a_j + n_1 b_1 + \cdots n_k b_k$$
with all coefficients positive integers and $\sum m'_i + \sum n_i < \sum m_i$. Furthermore, if this happens we can suppose that each of these coeficients is divisible by $f$. By subtracting the smaller of $m_i$ and $m'_i$ from the larger, we can assume that no $a_i$ appears on both sides of this equation. We take $T = \{b_1, \ldots, b_k\}$ and partition of $S$ into $S_1$ (those appearing on the left side of the equation) and $S_2$ (those not appearing on the left side of the equation). We take $T'$, $S'_1$ and $S'_2$ to be the corresponding sets in $W'$. It now transpires that we have written an element in the positive linear span of $S'_1$ as a positive linear combination of $S'_2 \cup T'$ and have a smaller coeficient sum using $S'_2 \cup T'$. But this happens exactly when $S'_1$ fails to lie on a face of $C$. This proves the second part of the proposition, and the third follows immediately. □

This leads to a few observations. Since the top dimensional faces are co-dimension 1, the union of the monoids,
$$\bigcup_{S \text{ is compatible}} M(S)$$
includes a finite index subgroup of $A$.

A "generic" generating set is one for which the points of $V'$ are in general position (modulo $a \mapsto a^{-1}$). In this case the faces of $C$ are simplices and there is no point of $V'$ in the interior of a face. Hence, when $S$ is compatible $\pi$ is an isomorphism so that $p_\mathcal{A}|_{M(S)} = p_{\mathbb{Z}^n_{\geq 0}}$, and when $S$ is maximal compatible, $n$ is the rank of $A$. Thus, in the generic case, there are pieces of $G$ so that for any piece a finite index subgroup of $G$ meets that piece in a simplicial cone and the Pascal's function there looks exactly like the standard one.

## Hyperbolic groups

There is a general method for finding the Pascal's function of a word hyperbolic group.

**Theorem.** *Let $G$ be a word hyperbolic group and let $\mathcal{G} = \{g_1, \ldots g_k\}$ be a generating set for $G$. Then there are $m \times m$ matrices, $M_1, \ldots, M_k$ and vectors*

$u = [u_1 \ldots u_m]$ and $v = \begin{bmatrix} v_1 \\ \vdots \\ v_m \end{bmatrix}$ with the following property: If $g \in G$ and $g_{i_1} \ldots g_{i_n}$ is any geodesic for $g$, then

$$p_{\mathcal{G}}(g) = u M_{i_1} \ldots M_{i_n} v.$$

*Proof.* Let $L$ be the set of all $\mathcal{G}$ geodesics. $L$ is the language of an automatic structure, and in particular, there is a finite state automaton $F$ which determines whether two geodesics represent the same element of $G$. This finite state automaton can be seen as a finite labelled directed graph: the vertices of the graph correspond to the states of the machine, the edges correspond to the transitions and the labels on those edges correspond to the input letters mediating those transitions. Each letter here is a pair $(a, a')$ where each of $a$ and $a'$ is either blank or an element of $\mathcal{G}$, and they are not both blank. (In fact, in this machine no edge leading to an accept state has a blank for either $a$ or $a'$.)

If we now fix an element $g$ and a geodesic $g_{i_1} \ldots g_{i_n}$ for $g$, then $p_{\mathcal{G}}(g)$ is the number of words $w = a_1 \ldots a_n$ so that the pairs $(g_{i_1}, a_1) \ldots (g_{i_n}, a_n)$ label a path starting from the start state of $F$ to an accept state of $F$. To count these we do the following. We take $m$ to be the number of states of $F$ and suppose that these are enumerated $s_1, \ldots, s_m$. (We assume $s_1$ is the start state.) We define $M_i$ to be the $m \times m$ matrix so that $m_{ij}$ gives the number of edges from state $i$ to state $j$ bearing a label of the form $(g_i, a')$. We take $u = [1\ 0\ \ldots\ ]$. We take $v$ so that the $i^{\text{th}}$ entry is 1 if $s_i$ is an accept state and 0 otherwise.

A standard induction shows that this does what is required. That is, if $0 \leq r \leq n = \ell(g)$, we let $N = N_r = M_{i_1} \ldots M_{i_r}$. ($N_0 = I$.) Then $n_{ij}$ gives the number of paths from $s_i$ to $s_j$ labelled by words of the form $(g_{i_1}, a_1) \ldots (g_{i_r}, a_r)$. Pre- and post-multiplication by $u$ and $v$ sum over paths from the start state to accept states. We leave the details to the reader. □

Bartholdi [**B**] has similar and more efficient methods in the case where $G$ is a hyperbolic surface group.

**Proposition.** *Suppose that $G = A * B$ and that $\mathcal{G} = \mathcal{A} \cup \mathcal{B}$, where $\mathcal{A}$ and $\mathcal{B}$ are generating sets for $A$ and $B$ respectively. Then*

$$p_{\mathcal{G}}(a_1 b_1 \ldots a_k b_k) = p_{\mathcal{A}}(a_1) p_{\mathcal{B}}(b_1) \ldots p_{\mathcal{A}}(a_k) p_{\mathcal{B}}(b_k).$$

*Proof.* A $\mathcal{G}$ geodesic consists of $\mathcal{A}$ and $\mathcal{B}$ geodesics for its factors. □

If the Pascal's function of a group graph is identically 1, then there is a unique geodesic to each group element. It is easy to arrange for this to happen in any finite group: we take the entire group as the generating set. Likewise this happens in a free group if we take our generating set to be a basis. It now follows that an arbitrary product of free and finite groups has a generating set in which the Pascal's function is identically 1. This raises the following

**Question.** *Suppose $G$ has a generating set for which Pascal's function is identically 1. Does it follow that $G$ is a free product of free groups and finite groups?*

Papasoglu has given a partial answer to this in [**P**] where he has shown that if a group is hyperbolic and has Pascal's function identically 1 then it is virtually free.

We prove the following.

**Theorem.** *Suppose $G = \langle \mathcal{G} \rangle$ is virtually infinite cyclic and that $p_{\mathcal{G}}$ is identically 1. Then either $G = \mathbb{Z}$ and $\mathcal{G}$ is a single generator, or $G$ is the infinite dihedral group $\mathbb{Z}_2 * \mathbb{Z}_2$ and $\mathcal{G}$ consists of two involutions.*

As we will observe below, it is sufficient to assume that $p_{\mathcal{G}}$ is bounded.

*Proof.* Since $G$ is virtually cyclic, it is word hyperbolic. Hence there is a finite state automaton $F$ whose language is the entire language of geodesics. Since the set of geodesics is infix closed, we can assume that every vertex is both a start state and an accept state. Since $G$ is infinite, $F$ has a loop. Assume the label on this loop is $y = y_1 \ldots y_j$. For some power $m$, $\langle \overline{y}^m \rangle$ is normal in $G$. We consider the word $y^m y_1$. The word $y^m y_1$ is geodesic, and we either have $\overline{y^m y_1} = \overline{y_1 y^m}$ or $\overline{y^m y_1} = \overline{y_1 y^{-m}}$.

*Case 1:* $\overline{y^m y_1} = \overline{y_1 y^m}$. Then $y_1 y^m$ is also geodesic, and necessarily equal $y^m y_1$. This implies that $y$ is a power of $y_1$. In this case we will call $y_1 = t$ and $y = t^n$ for some $n$.

*Case 2:* $\overline{y^m y_1} = \overline{y_1 y^{-m}}$. Now $y^m y_1$ is a geodesic, and this evaluates to the same element as $y_1 y^{-m}$. Since these both have the same length, the latter is also geodesic, so these two words are necessarily equal. But $y^m y_1$ ends in $y_1$ and $y_1 y^{-m}$ ends in $y_1^{-1}$. Evidently $y_1 = y_1^{-1}$ and $\overline{y_1} = y_1^{-1}$.

Now $j \geq 2$, since $y^2$ is geodesic, while $\overline{y_1^2} = 1$. We look at the loop labeled $y' = y_2 \ldots y_j y_1$ based at the next vertex of the loop $y$. Since $y'$ is a cyclic conjugate of $y$, and $\langle \overline{y}^m \rangle$ is normal, $\langle \overline{y}^m \rangle = \langle \overline{y'}^m \rangle$. In particular, $\langle \overline{y'}^m \rangle$ is normal. Performing the same argument as before, we have $\overline{y'^m y_2} = \overline{y_2 y'^{e_2 m}}$, where $e_2 = \pm 1$. But $e_2 = 1$ is impossible, since we then have (as in case 1) $y'$ is a power of $y_2$ whence $y'$ is a power of $y_1$ which is of order 2. Consequently $\overline{y'^m y_2} = \overline{y_2 y'^{-m}}$, and, as in case 2, $y_2$ has order 2.

But now we observe that $\overline{y^m y_1 y_2} = \overline{y_1 y_2 y^m}$ so that $y^m y_1 y_2 = y_1 y_2 y^m$, and thus $y$ is a power of $y_1 y_2$. In this case we call $y_1$ and $y_2$ $r$ and $s$ respectively and have $y = (rs)^n$ for some $n$.

Let the loops of $F$ bear the labels $v_1, \ldots, v_q$. Then all of $\overline{v_1}, \ldots, \overline{v_q}$ have powers lying in a common normal subgroup $Z = \langle z \rangle < G$. Thus, for each $i$ some power of $\overline{v_i}$ is a power of $z$ or of $z^{-1}$. This divides the set of loops of $F$ into two equivalence classes. Now if $v_i$ and $v_j$ fall in the same equivalence class, $\overline{v_i}$ and $\overline{v_j}$ have a common power, so $v_i$ and $v_j$ are themselves powers of common word. In particular they are both labelled by a power of either $t$, $t^{-1}$, $rs$ or $sr$. Furthermore, if $v_i$ labels a loop, so does its inverse, since the inverse of a geodesic is also a geodesic. Thus all the loops of $F$ are labelled by positive and negative powers of $t$ or all the loops of $F$ are labelled by powers of $rs$ and its inverse $sr$.

Now the set of all geodesics is a regular language in which the number of words of length $n$ is bounded by a linear function of $n$. It is a standard result that such

languages are finite unions of the form

$$\cup_i \{u_i v_i^m w_i \mid m \geq 0\}.$$

To finish the proof it only remains to see that if each $v_i$ is a power of $t$, then so are each $u_i$ and $v_i$ and that if each $v_i$ is a power of $rs$, then each $u_i$ and $v_i$ consists only of $r$'s and $s$'s. (This is certainly true for any $u_i$ or $v_i$ that is empty!)

*Case 1*: Each $v_i$ is a power of $t$. We repeat the arguement of case 1 above using using the last letter of $u_i$ or the first letter of $w_i$ in the rôle of $y_1$. We then repeat this peeling off successive letters of $u_i$ and $w_i$, thus showing that each of these consists only of $t^{\pm 1}$'s.

*Case 2*: Each $v_i$ is a power of $rs$. We suppose $v_i$ is a positive power of $rs$. Then the last letter of $u_i$ conjugates a power of $v_i$ to its inverse and is thus $s$. The last two letters of $u_i$ (if there are two) conjugate a power of $rs$ to itself, and are thus $rs$. Thus each nonempty $u_i$ is an alternating word in $r$ and $s$ ending in $s$ and likewise each $w_i$ is an alternating word in $r$ and $s$ beginning in $r$. □

We can weaken the supposition that $p_{\mathcal{G}} = 1$ to the supposition that $p_{\mathcal{G}}$ is bounded. For if we can move (say) $x$ through $y^m$ (where $y$ labels a loop) either preserving or reversing sign, but giving a different word, then we can change $x(y^m)^k$ into any of $k$ different geodesic words for the same element.

Department of Mathematics
University of Melbourne
Parkville, VIC 3052
Australia